\documentclass{amsart}

\usepackage{amsmath,amssymb,enumerate,amsthm}
\usepackage{epsfig}

\swapnumbers
\newtheorem{proposition}{Proposition}
\newtheorem{theorem}[proposition]{Theorem}

\newtheorem{lemma}[proposition]{Lemma}
\newtheorem{corollary}[proposition]{Corollary}

\newtheorem{definition}[proposition]{Definition}

\numberwithin{proposition}{section}
\numberwithin{equation}{section}

\newcommand{\occult}[1]{}

\newcommand\eps{\epsilon}

\newcommand\NN{{\mathbb N}}

\renewcommand\top{{\operatorname{top}}}

\begin{document}

\title[No Measure of Maximal Entropy]{A continuous, piecewise affine surface map with\\
no measure of maximal entropy}

\begin{abstract}
It is known that piecewise affine surface homeomorphisms always
have measures of maximal entropy. This is easily seen to fail
in the discontinuous case. Here we describe a piecewise affine, globally
continuous surface map with no measure of maximal entropy.
\end{abstract}

\author{J\'er\^ome Buzzi}

\address{Jerome Buzzi\\C.N.R.S. \& D\'epartement de Math\'ematique\\
Universit\'e Paris-Sud\\91405 Orsay Cedex\\France}

\email{jerome.buzzi@math.u-psud.fr}

\urladdr{www.jeromebuzzi.com}


\maketitle

\section{Introduction}

\subsection{Measures of Maximal Entropy}

The complexity of the orbit structure of a dynamical system is
reflected by its \emph{topological entropy}. Entropy can also be defined at
the level of invariant probability measures. These two levels
are related by the variational principle, i.e., for any continuous
self-map of a compact metric space:
 $$
     h_\top(f) = \sup_\mu h(f,\mu).
 $$
This brings to the fore \emph{maximal entropy measures}, i.e., those having "full
complexity" in the following sense: 
 $$
   h(f,\mu)=h_\top(f).
 $$
Such measures may fail to exist, e.g., for any $r<\infty$, there are $C^r$ smooth interval
maps with non-zero topological entropy and no maximal entropy
measures \cite{BuzziSIM,RuetteEx}. However, building on Yomdin's theory
\cite{Yomdin} of smooth mappings, Newhouse has shown the following

\begin{theorem}[Newhouse \cite{Newhouse}]
If $f$ is $C^\infty$ self-map of a compact manifold, then
$\mu\mapsto h(f,\mu)$ is upper semicontinuous as a function on
the compact set of invariant probability measures endowed with
the weak star topology. In particular, there exists a maximal
entropy measure.
\end{theorem}

\subsection{Piecewise Affine Transformations}
This does not apply to the following simple class of transformations,
even under the assumption of global continuity:

\begin{definition}
A map $T:M\to M$ is said to be {\bf piecewise affine} if
 \begin{enumerate}
 \item $M$ is admits an affine atlas (i.e., a set of charts whose change 
of coordinates are affine diffeomorphisms);
 \item there exists a finite partition of $M$ whose elements $A$ satisfy:
(i) $A$, resp. $TA$, are each contained in the domain of a chart $\chi$, resp. $\chi'$,
of the affine atlas; (ii) $\chi'\circ T\circ\chi^{-1}|A$ is the restriction of an
affine map to some open subset of some $\mathbb R^d$.
 \end{enumerate}
\end{definition}

However, Newhouse observed
\cite{NewhousePerso,BuzziPWAH} that the
above property nevertheless holds for \emph{piecewise affine surface
homeomorphisms}. This follows from the sub-exponential rate
at which discontinuities can accumulate when one iterates the
map: the \emph{multiplicity entropy} introduced in \cite{BuzziAffine}
is zero for such maps.

Additionally, the set of maximal entropy measures of such
transformations has been shown \cite{BuzziPWAH} to be a finite-dimensional
simplex whenever $h_\top(f)\ne0$.

It is easy to see that the {\bf finiteness property fails
for piecewise affine continuous maps}. Indeed, it is enough to
consider the direct product of the identity on some interval
with a piecewise affine, globally continuous interval map with
nonzero entropy.

\subsection{Main Result}
In this note we show that {\bf existence also fails} to
hold generally for such maps:

\begin{theorem}
There exists a piecewise affine, globally continuous map $T$ of $[0,1]^2$ satisfying
 $
    h_\top(T)=\sup_\mu h(T,\mu)=\log2
 $
with no maximal entropy measure, i.e., no measure achieving this supremum.
\end{theorem}

This map will be explicitly described. Taking the direct product of the
above map with the identity of a cube of the proper dimension, one immediately
obtains the following:

\begin{corollary}
For any integer $d\geq2$, there exists a piecewise affine, globally continuous map $T$ of $[0,1]^d$ satisfying
 $
    h_\top(T)=\sup_\mu h(T,\mu)=\log2
 $
with no maximal entropy measure, i.e., no measure achieving this supremum.
\end{corollary}

\subsection{Comments}
We can compare globally continuous, piecewise affine maps to related classes
for which this problem has been studied. 

On the one hand, existence is known to fail if one removes:
 \begin{itemize}
  \item either the continuity assumption. We refer to
\cite{BuzziAffine,KruglikovRypdal1,KruglikovRypdal2} for examples and more
discussion, including failure of the variational principle.
  \item or the affine assumption. It is easy to construct a globally continuous,
  piecewise quadratic map without a measure of maximal entropy \cite{BuzziPWAH}. 
  This construction is derived from examples of $C^\infty$ maps at which the topological entropy
  fails to be lower semi-continuous.
 \end{itemize}
On the other hand, a classical theorem of Hofbauer \cite{Hofbauer} asserts 
that existence and finiteness hold for piecewise affine interval maps 
(in fact, piecewise monotone maps) with non zero entropy. 

\section{General Description of the Map}\label{sec:map-desc}

\subsection{Key Properties}

The map $T$ is defined on a parallelogram $Q:=ENWS$
(we write $X_1\dots X_k$ for the compact polygon with sides $[X_1X_2],\dots,[X_kX_{1}]$). 

Define the vertical cone
 $$
   \mathcal C^s:=\left\{\left(\begin{matrix}x\\y\end{matrix}\right):|x|\leq 2|y|\right\}.
 $$
We say that $\mathcal C^s$ is \emph{stable} for some map, if the differential at any point of
the inverse of that map (where this differential exists) sends $\mathcal C^s$ into itself.

\begin{figure}
\centering
\includegraphics[width=10cm]{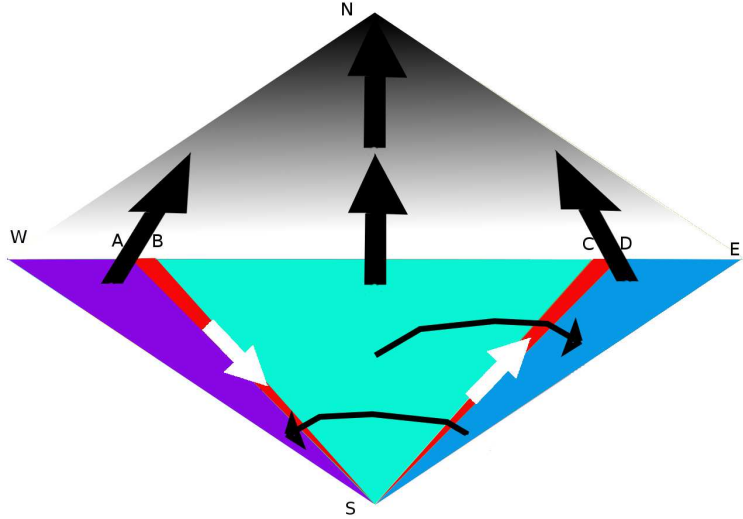}
\caption{Main zones for the dynamics of the example $T$}\label{fig:scheme}
\end{figure}

Let us describe the key properties of the map $T$ (see Fig. \ref{fig:scheme}):

\begin{enumerate}
 \item\label{key-NS} the poles $N$ and $S$ are fixed points;
 \item\label{key-half} all points in the top half $NWE$ (the grey triangle in Fig. \ref{fig:scheme})
 are attracted to $N$;
 \item\label{key-markov} each one of the red triangles $ABS$ and $CDS$ is mapped to the large
 triangle $ADS$;
 \item\label{key-ABS} $T|ABS$ divides the $y$-coordinate by a factor $\leq 2$;
 \item\label{key-CDS} $T|CDS$ multiplies the $y$-coordinate by a factor of $2$;
 \item\label{key-cone} $\mathcal C^s$ is stable for the map on $ABS\cup CDS$;
 \item\label{key-transverse} on $ABS\cup CDS$, the map preserves the horizontal and expands horizontal vectors by a factor $\geq 4$; 
 \item\label{key-folding} the middle, green-blue part $BCS$ is mapped into
  the blue, right side $DES$, except for a part which is mapped into $NEW$ after one or two
 iterations of $T$;
 \item\label{key-right} $DES$ is mapped into the purple left side $WAS$ and $NWE$;
 \item\label{key-left} all orbits starting in $WAS$ eventually converge to fixed points.
\end{enumerate}
To be precise, we must state the following corrections, involving the exact partition defined below:

Properties (\ref{key-markov})-(\ref{key-transverse}) do not hold on the top parts
$ABA^tB^t$ and $CDC^cD^c$ which are mapped into $NEW$. The contraction is exactly by a factor of $2$
on the lower part, $A^cB^cS$, of $ABS$.

\subsection{Precise definition of the map}

The above is realized as a piecewise affine, globally continuous map $T:Q\to Q$. The parallelogram $Q=NWES$ is partitioned into 26 triangles, on each of which the map is affine.

\begin{figure}
\centering
\hskip-1.5cm\includegraphics[width=14cm]{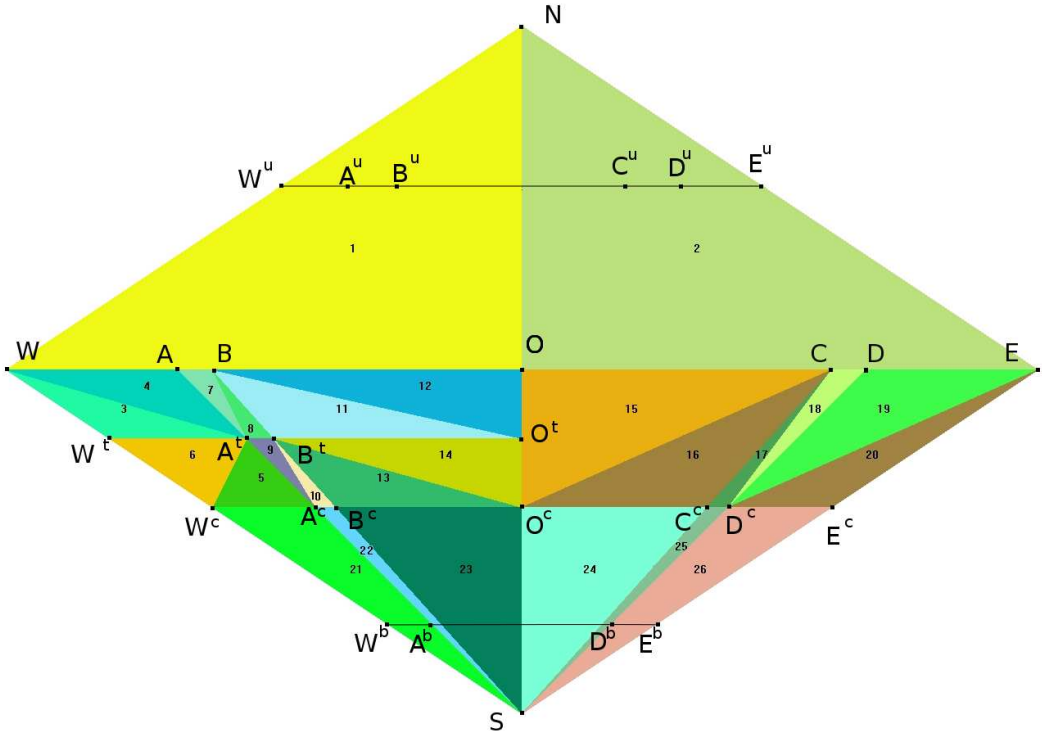}
\caption{Partition and special points for the map $T$.}\label{fig:partition}
\end{figure}

The $x$-coordinates of the vertices of this partition which lie on the 
diagonal $y=y^1$ are given in Table \ref{tab:key-vert}.
The top and bottom points of $Q$ are $N(0,2)$ and $S(0,0)$.  The other vertices of the
partition are obtained from those on $y=y^1$ by homotheties centered at $S$ or $N$ yielding
points on $y=y^u,y^t,y^c,y^b$ with values in Table \ref{tab:key-lines}. We denote by
$A^t$ the vertex obtained from $A$ on the line $y=y^t$, etc.
The resulting partition and vertices are depicted in Fig. \ref{fig:partition}.

\begin{table}
\begin{tabular}{|c|c|c|c|c|c|c|c|c|c|}
\hline
Name & W    & A  & B    & O & C   & D & E\\
\hline
$x$		& -1.5 & -1 & -0.9 & 0 & 0.9 & 1 & 1.5\\
\hline
\end{tabular}
\medbreak
\caption{$x$-coordinates of the vertices on the line $y=y^1$.}\label{tab:key-vert}
\end{table}

\begin{table}
\begin{tabular}{|*{10}{c|}}
\hline
Name            & $y^1$  & $y^t$ & $y^c$ & $y^b$ & $y^u$ \\
\hline
$y$		& 1    & 0.8 & 0.5 & 0.25 & 1.5\\
\hline
\end{tabular}
\medbreak
\caption{$y$-coordinates of the horizontal lines used to define $T$.}\label{tab:key-lines}
\end{table}

The piecewise affine map is finally defined by the images of the vertices
of the partition, as given in Table \ref{tab:map-points}. 

\begin{table}
\begin{tabular}{|*{10}{c|}}
\hline
On $y=y^1$ & $W$ & $A$ & $B$ & $O$ & $C$ & $D$ & $E$ \\
\hline
Image  & $W^u$ & $A^u$ & $D^u$ & $E^u$ & $D^u$ & $A^u$ & $W^u$\\
\hline
\hline
On $y=y^t$ & $W^t$ & $A^t$ & $B^t$ & $O^t$ & --- & --- & --- \\
\hline
Image  & $W$ & $A$ & $D$ & $E$ & ---& --- & --- \\
\hline
\hline
On $y=y^c$ &  $W^c$ & $A^c$ & $B^c$ & $O^c$ & $C^c$ & $D^c$ & $E^c$ \\
\hline
Image  & $W^c$ & $A^b$ & $D^b$ & $E^c$ & $D$ & $A$ & $W$\\
\hline
\end{tabular}
\medbreak
\caption{Mapping of the points. Note: the remaining vertices, $S$ and $N$, are fixed points.}\label{tab:map-points}
\end{table}

\subsection{Proof of the Key Properties}\label{sec:map-check}

Key properties (\ref{key-NS}) and (\ref{key-half}) are immediate from
$T(S)=S$, $T(N)=N$ and the fact that $T|NWO$ and $T|NEO$ are affine
map with $W,O,E$ mapped strictly above the line $(WE)$. 

\begin{figure}
\centering
\includegraphics[width=10cm]{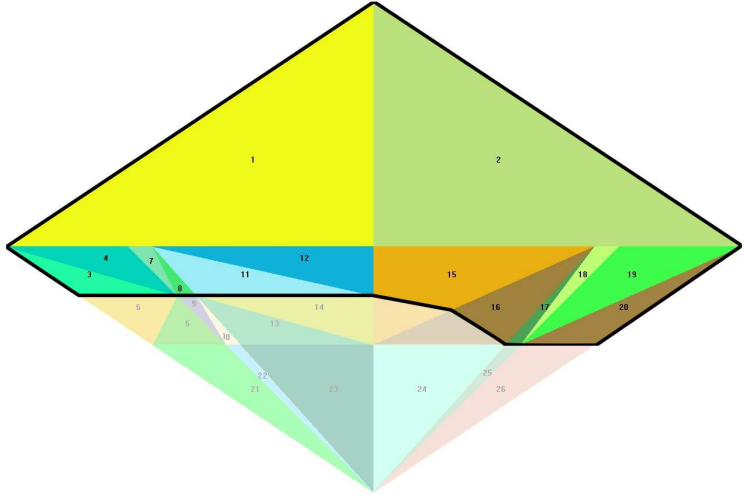}
\caption{$T^{-1}(NEW)$ delineated by a bold line.}\label{fig:preimage-top}
\end{figure}

A direct computation shows that the preimage of $NWE$ is the union of $NWE$ with $WW^tOO^t$,
$C^cE^cCE$ and some upper part of $OO^cCC^c$ ---see Fig. \ref{fig:preimage-top}.
Using this, the following shows the key properties (\ref{key-markov})-(\ref{key-transverse}).

\begin{figure}
\centering
\includegraphics[width=6cm]{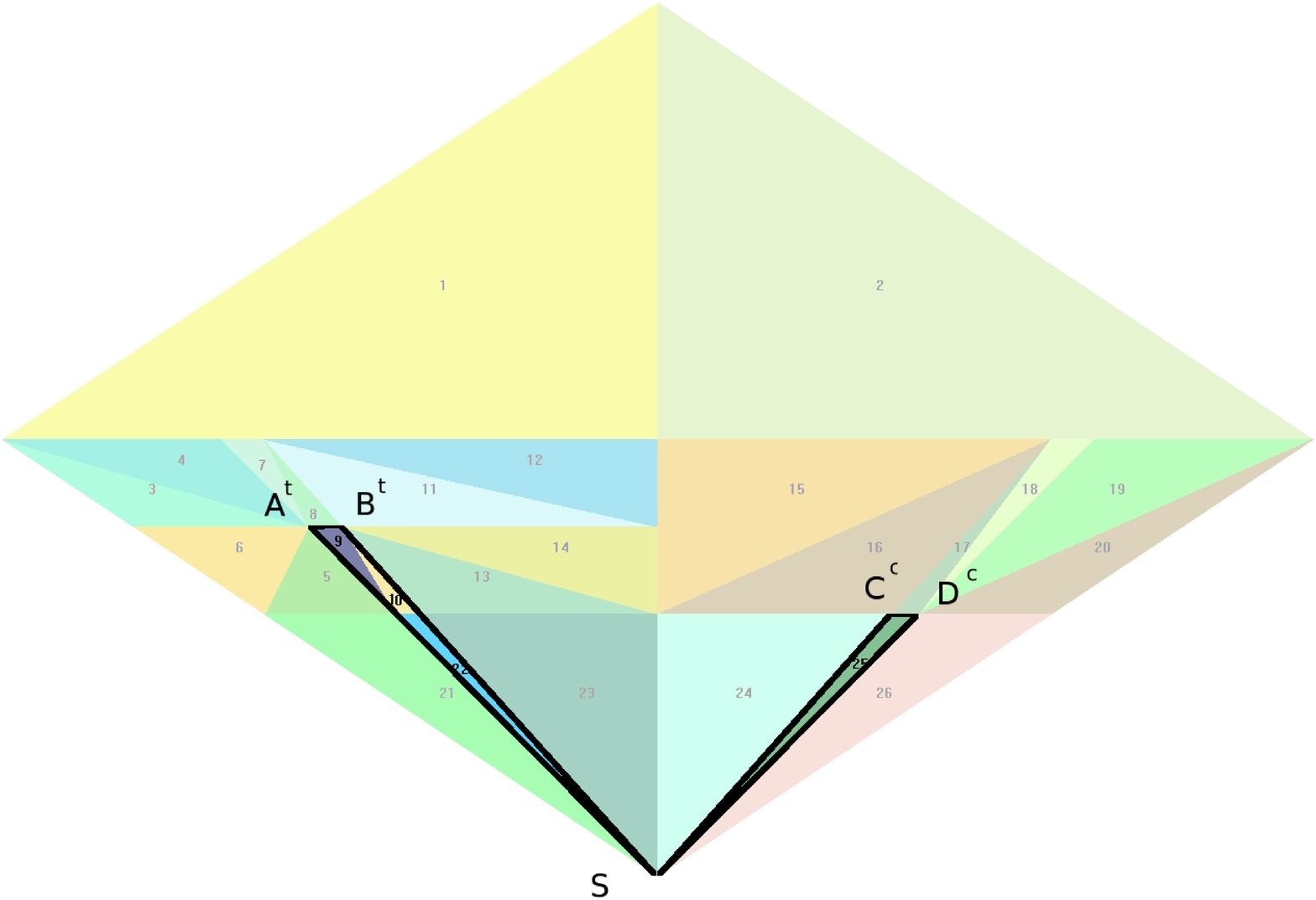}
\includegraphics[width=6cm]{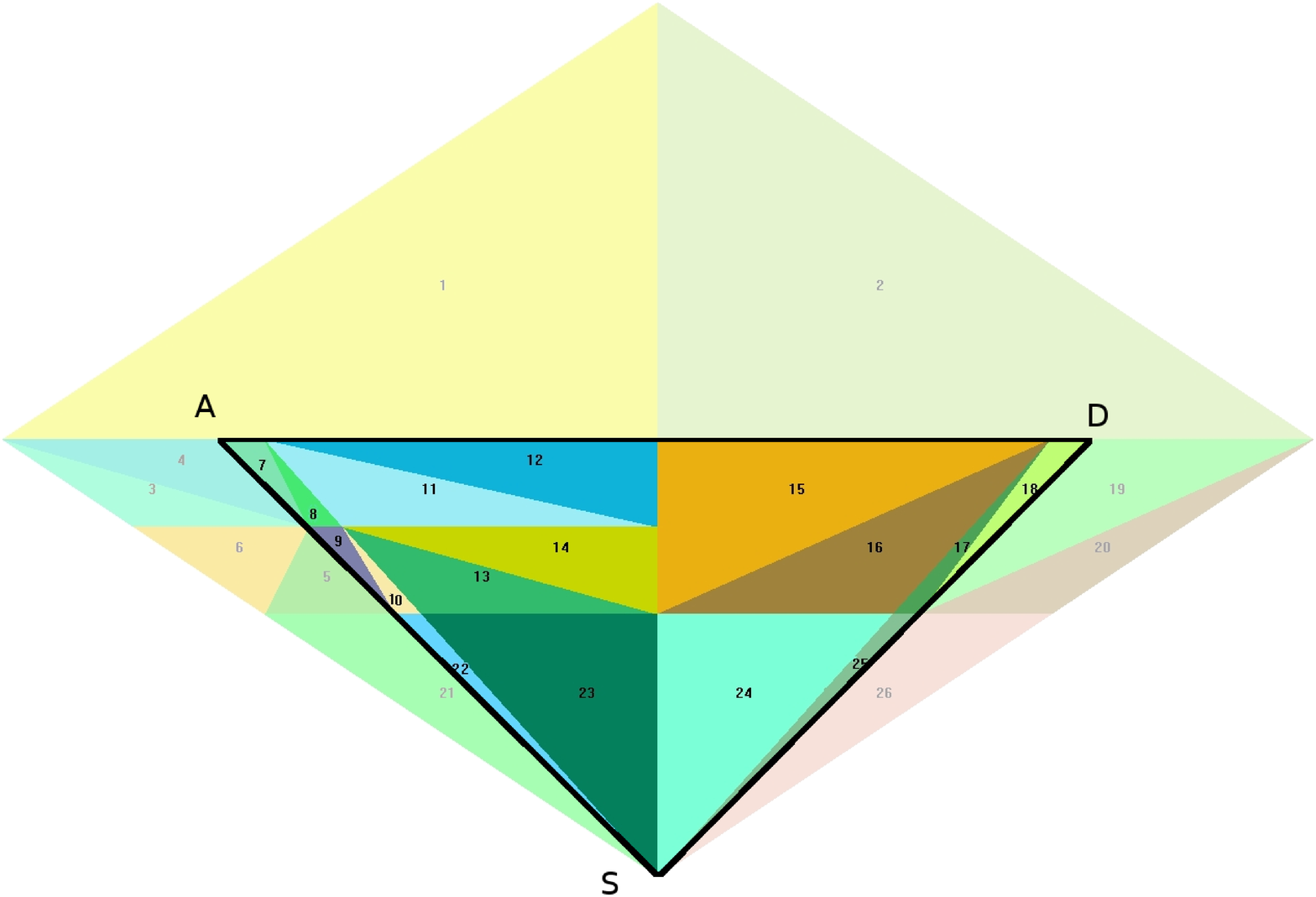}
\caption{On the left: the triangles $A^tB^tS$ and $C^cD^cS$. On the right: their (identical) image. The bold lines on the right are the images of those on the left.}\label{fig:strips}
\end{figure}

Key property (\ref{key-markov}) follows by inspection, see Fig. \ref{fig:strips}.
Looking at Table \ref{tab:mat} we see that $T'|A^tB^tS$ multiplies the $y$-coordinate
by $0.5$ or $2.5$, proving (\ref{key-ABS}), $0.5$ on $A^cB^cS$. (\ref{key-CDS})
is checked in the same way.

We are going to prove that $\mathcal C^s$ is stable for $T|A^tB^tS$ and $T|C^cD^cS$.
We have to check this property for the $4$ affine maps involved using the following:

\begin{lemma}\label{lem:cone-linear}
Let $A=\left(\begin{matrix} a & b\\ c& d\end{matrix}\right)$. The following is a sufficient condition
for the invariance $A^{-1}(K_C)\subset K_C$ where $K_C:=\{\left(\begin{matrix} x\\
y\end{matrix}\right):|x|\leq C|y|\}$:
 $$
    \gamma_1:=C|c/a|<1 \text{ and } \gamma_2:=\frac{|d|+|b|C^{-1}}{|a|-C|c|} \leq 1.
 $$
\end{lemma}

\begin{proof}
Let $\left(\begin{matrix} \xi'\\ \eta'\end{matrix}\right)=A\left(\begin{matrix} \xi\\ \eta \end{matrix}\right)$ with $|\xi'|\leq C|\eta'|$. We have
 $$
     \xi=a^{-1}\xi'-a^{-1}b\eta \text{ and }\eta'=c\xi+d\eta.
 $$
Hence, abreviating $|a|,|b|,\dots$ to $a,b,\dots$,
 $$
    |\xi| \leq a^{-1}|\xi'|+a^{-1}b|\eta| \leq C a^{-1}(c|\xi|+d|\eta|)+a^{-1}b|\eta|.
 $$
Thus,
 $$
     (1-a^{-1}cC)|\xi| \leq (a^{-1}dC+a^{-1}b) |\eta|.
 $$
If the first factor is positive this is equivalent to
 $$
    |\xi|\leq \frac{a^{-1}dC+a^{-1}b}{1-a^{-1}cC}.
 $$
The Lemma follows.
\end{proof}

\begin{table}
\begin{tabular}{|*{5}{c|}}
\hline
Triangle  & $A^cB^cS$ &  $A^tB^tA^c$        & $A^cB^tB^c$ & $C^cD^cS$ \\
\hline 
Matrix    & $10\quad 9.5$ & $25\quad 22.5$  & $10\quad 11.5$ & $-40\quad 38$\\
          & $\;0\quad 0.5$ & $\;0\quad 2.5$ & $\;0\quad 2.5$ & $0\quad 2$\\
\hline
$\gamma_1$& $0$  & $0$ & $0$ & $0$ \\
\hline
$\gamma_2$& $0.525$ & $0.55$ & $0.825$& $0.525$\\
\hline
\end{tabular}
\medbreak
\caption{Stability of the cone.}\label{tab:mat}
\end{table}

The matrices of the
linear parts of the $4$ affine maps are listed in Table \ref{tab:mat} together with
the quantities denoted above $\gamma_1,\gamma_2$ for each of those.
Thus (\ref{key-cone}) holds.

Property (\ref{key-transverse}) also follows immediately from the matrices in Table \ref{tab:mat}
(the left lower entry is zero and the left upper entry is bigger than $4$ in absolute
value).

\begin{figure}
\centering
\includegraphics[width=6cm]{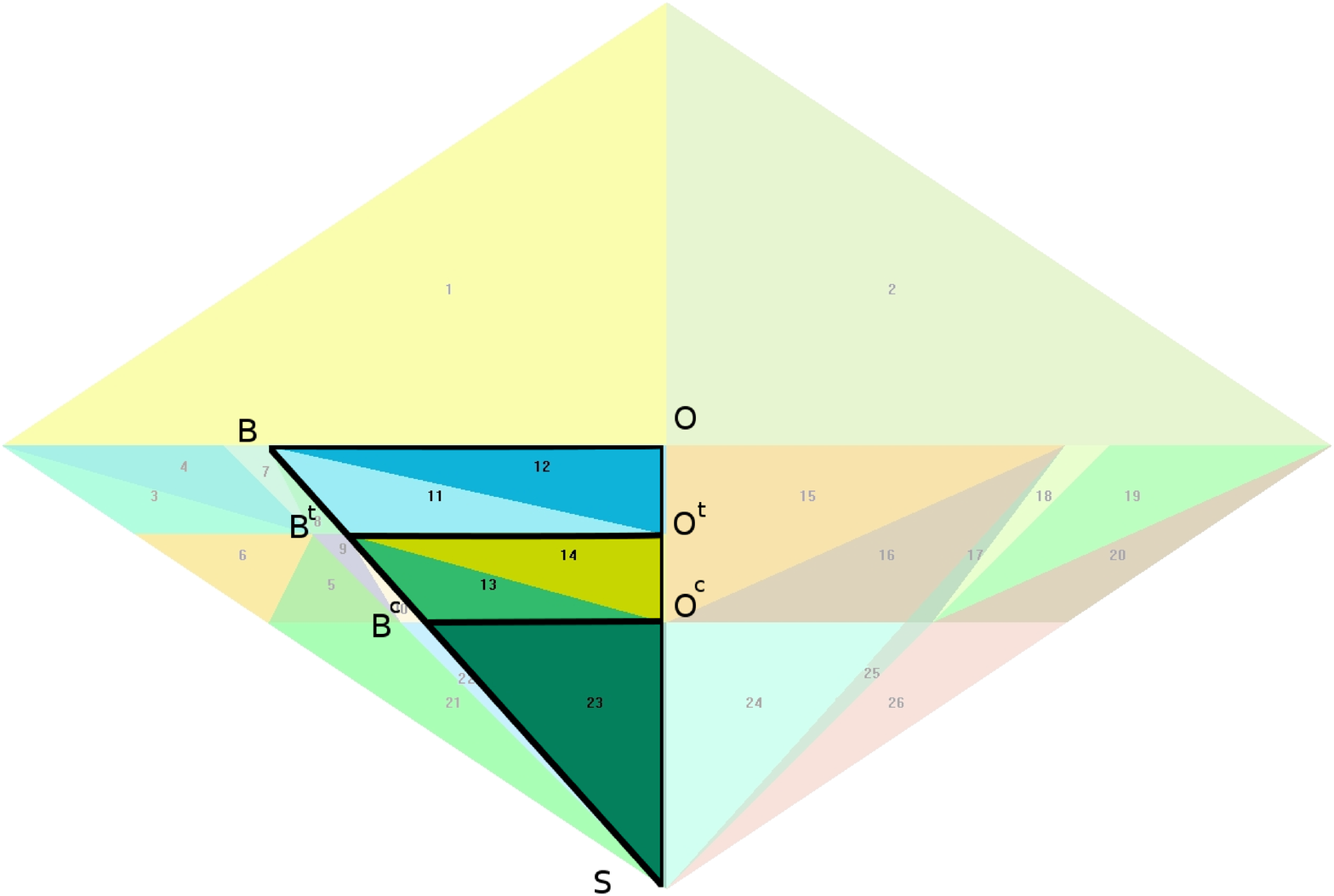}
\includegraphics[width=6cm]{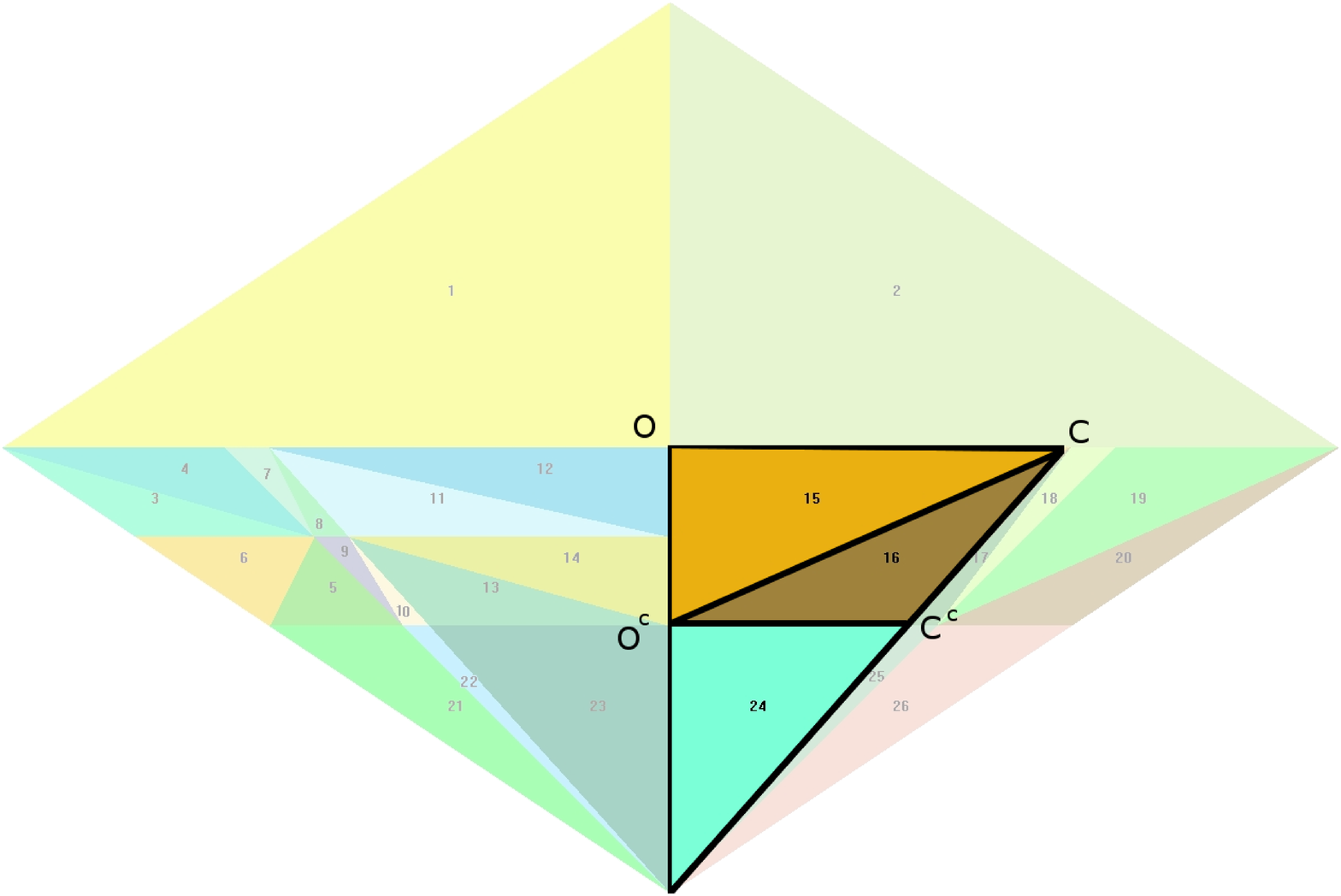}
\caption{The two halves of the folding zone: $BOS$ on the left and $OSC$ on the right.}\label{fig:middle}
\end{figure}

\begin{figure}
\centering
\includegraphics[width=6cm]{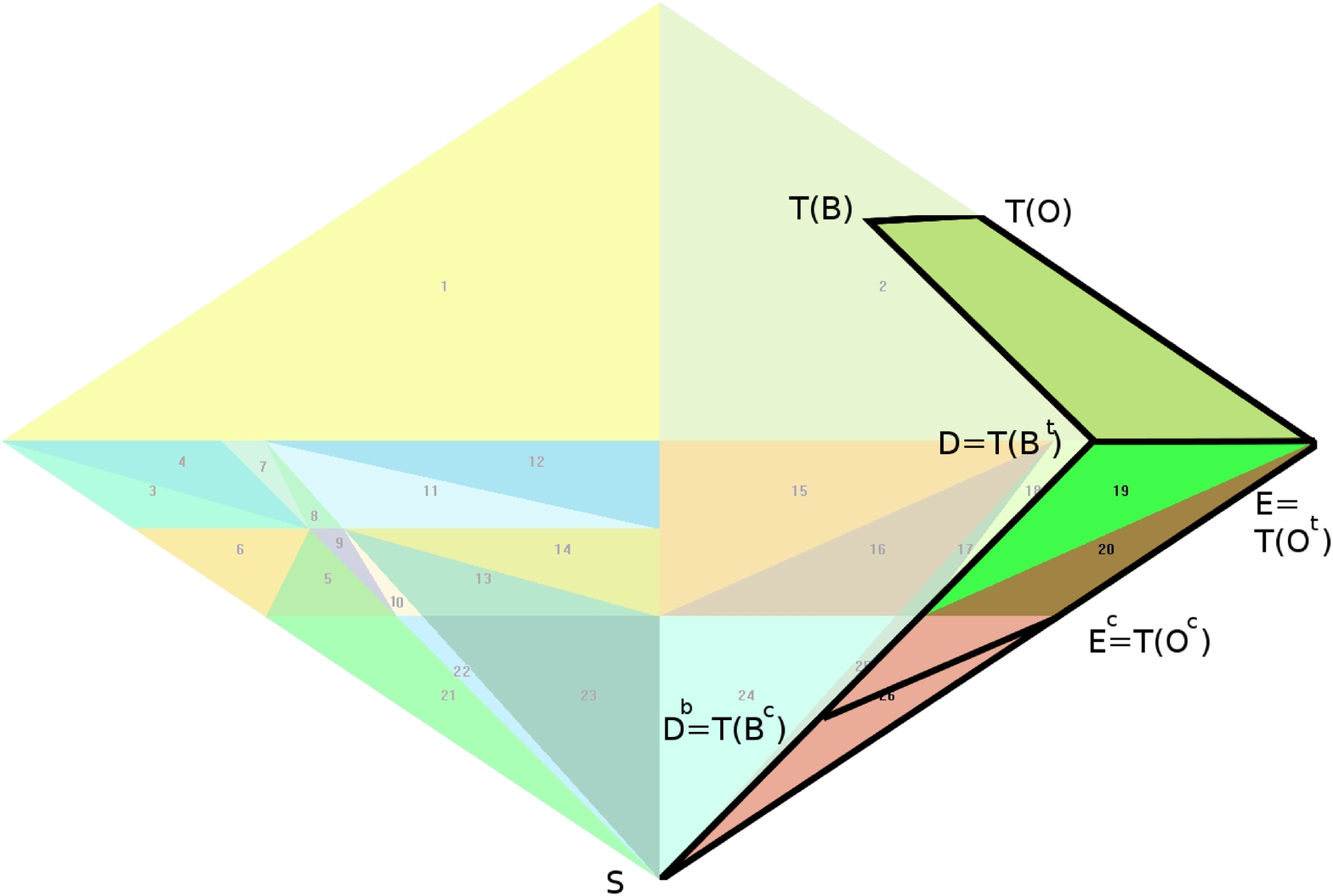}
\includegraphics[width=6cm]{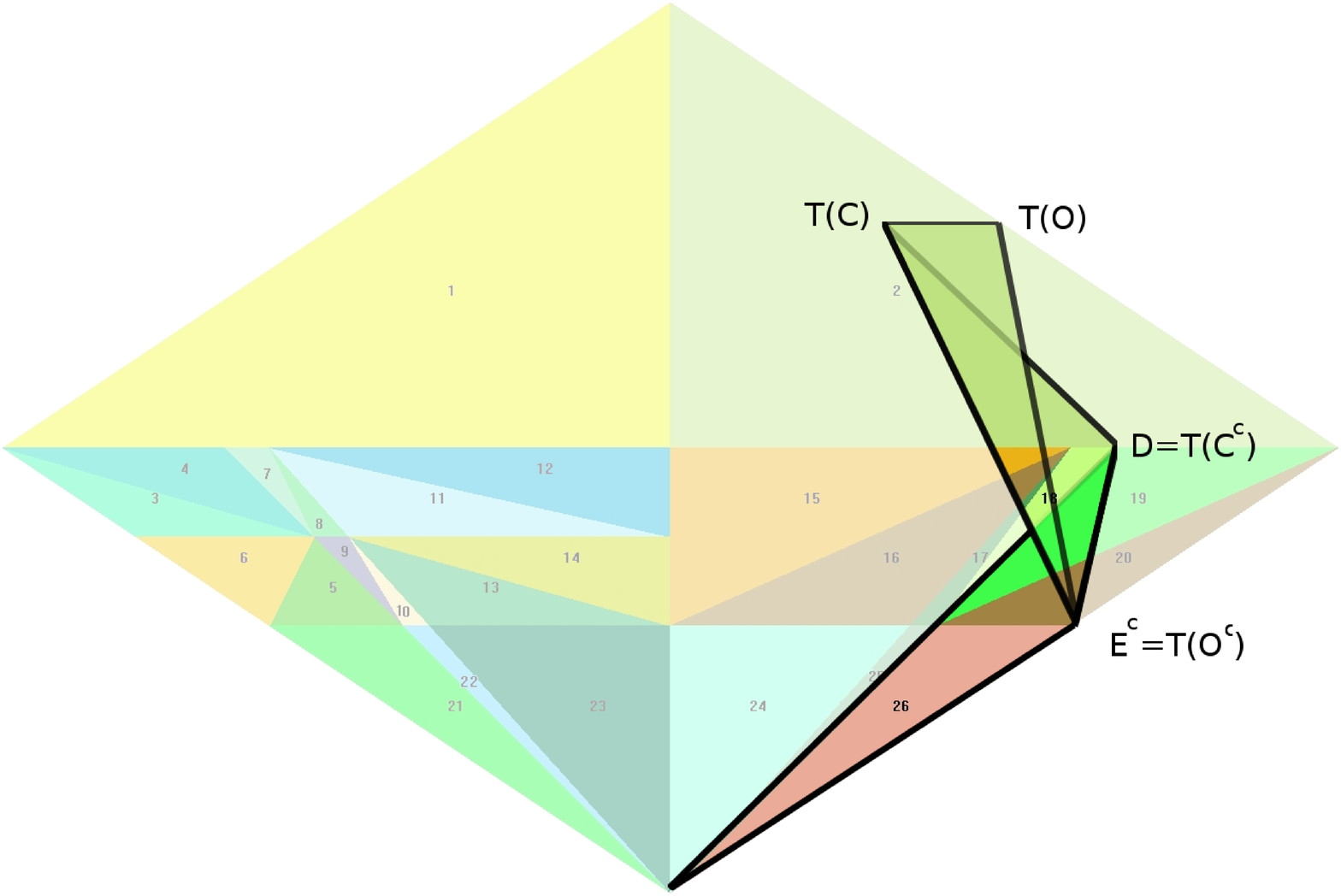}
\caption{The images of two halves of the folding zone: $T(BOS)$ on the left and $T(OSC)$ on the right. The bold lines are the images of those in Fig. \ref{fig:middle}.}\label{fig:middle-image}
\end{figure}

To prove that the "folding zone" $BCS$ is mapped to the right of $CDS$ except for
the part that ends up in $NEW$ in one or two iterations, we decompose it into its left half $BOS$ and right half
$OSC$ ---see Fig. \ref{fig:middle}. The images are given in Fig. \ref{fig:middle-image}.
We see that $T(BOS)\subset DES\cup NEW$ and the image of
$T(OSC)\subset DES \cup NEW\cup\Delta$ where $\Delta\subset OO^cCC^c$ is itself mapped
into $NEW$ according to Fig. \ref{fig:preimage-top}. This proves property (\ref{key-folding}).

\begin{figure}
\centering
\includegraphics[width=6cm]{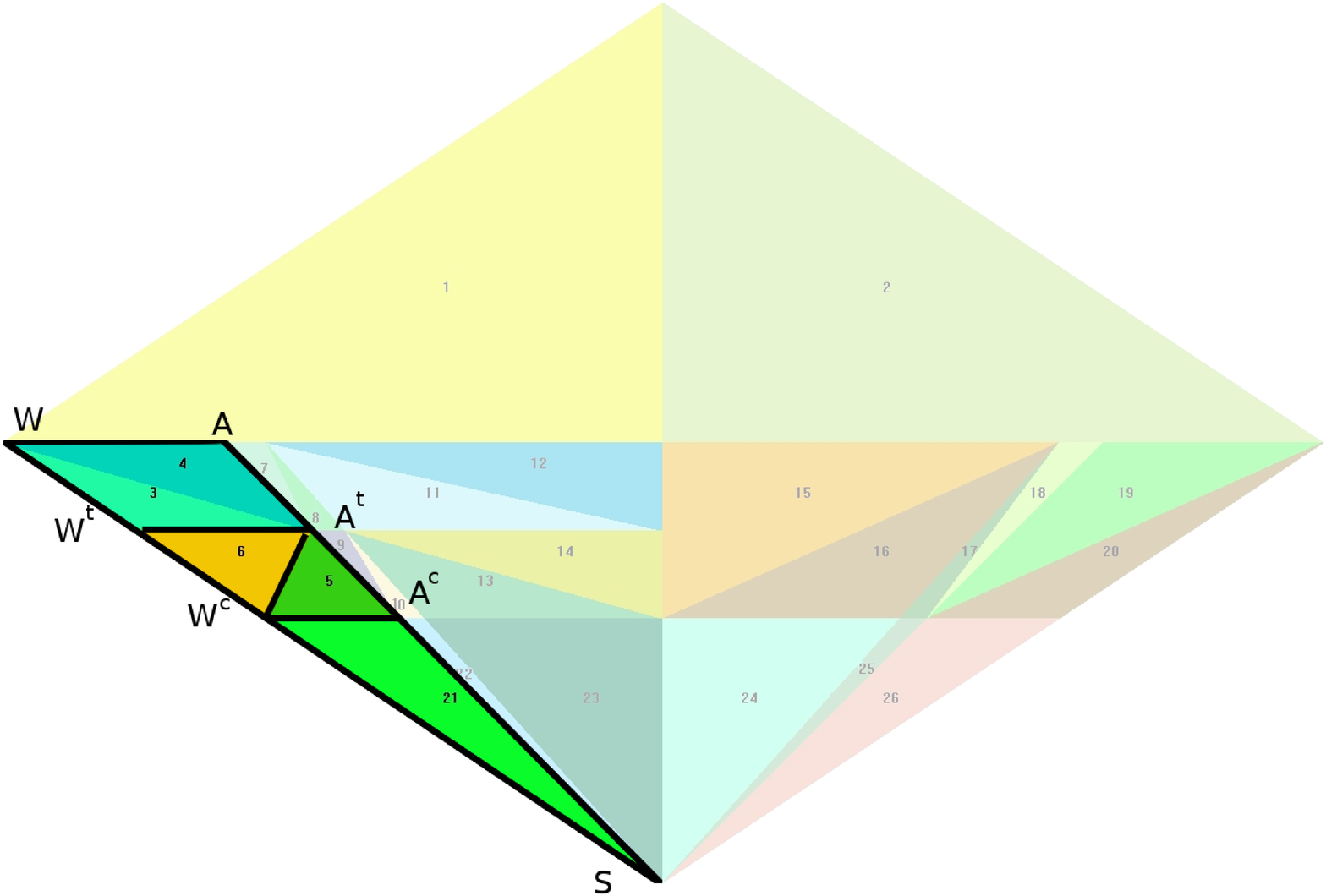}
\includegraphics[width=6cm]{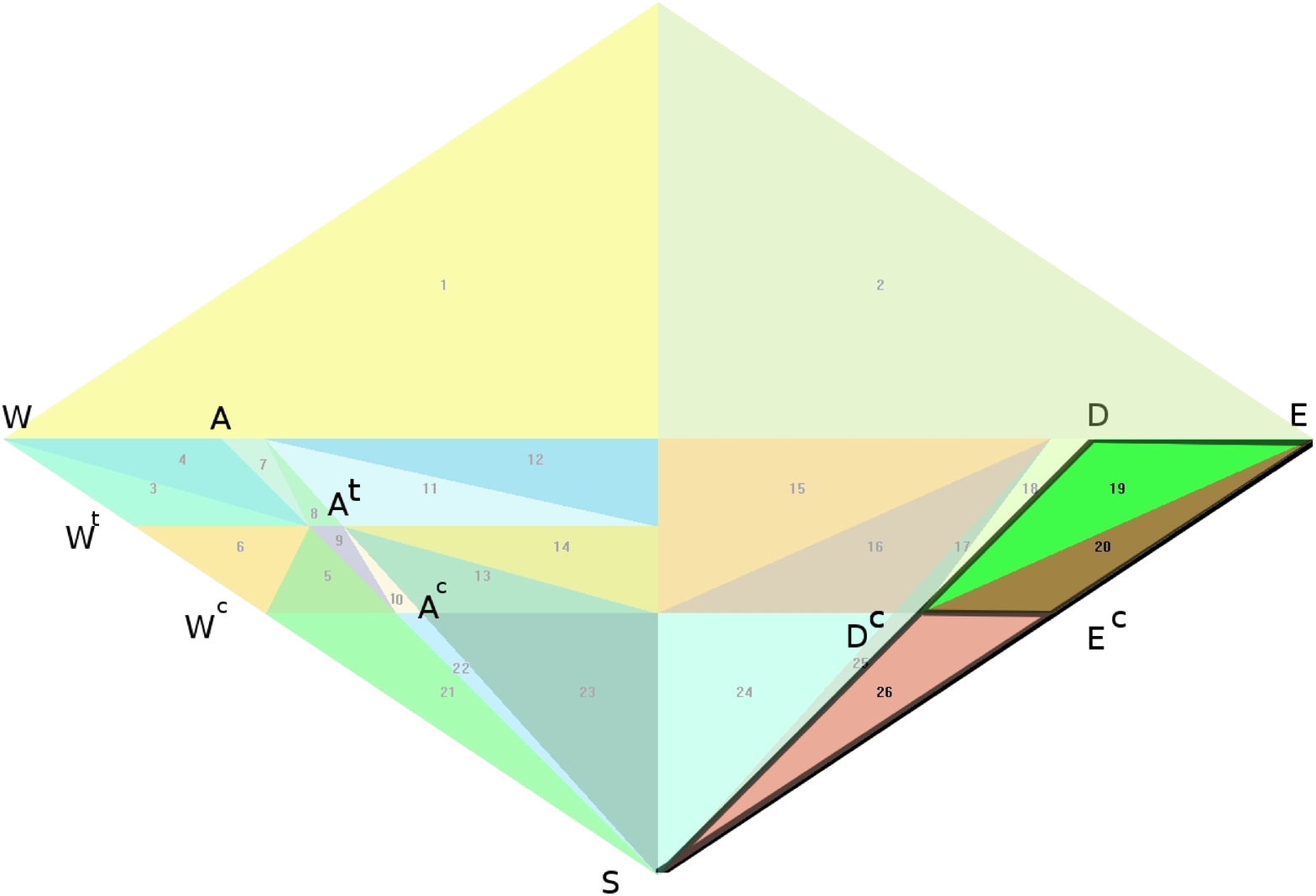}
\caption{The two triangles $(WAS)$ on the left and $(DES)$ on the right.}\label{fig:left-right}
\end{figure}

\begin{figure}
\centering
\includegraphics[width=6cm]{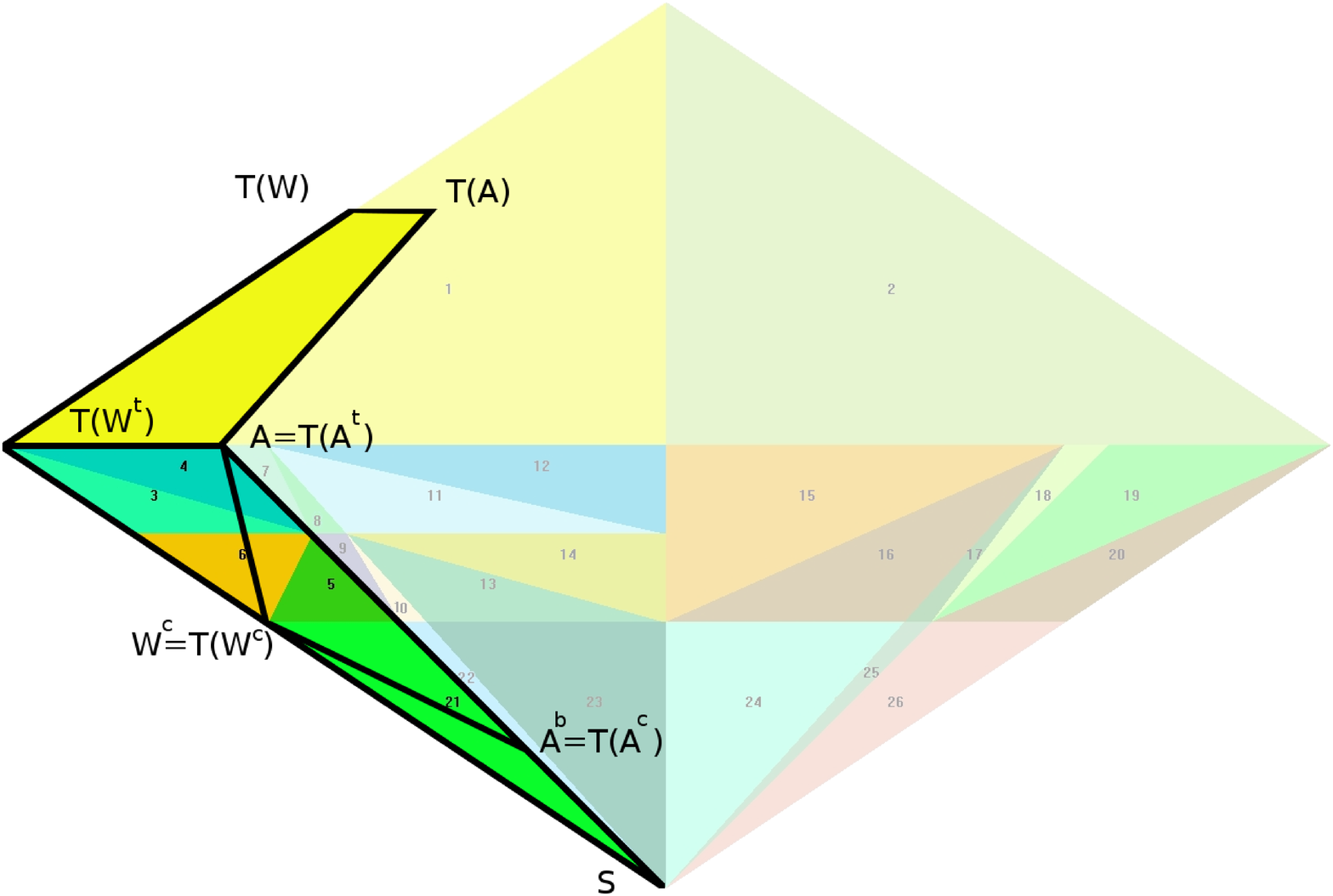}
\includegraphics[width=6cm]{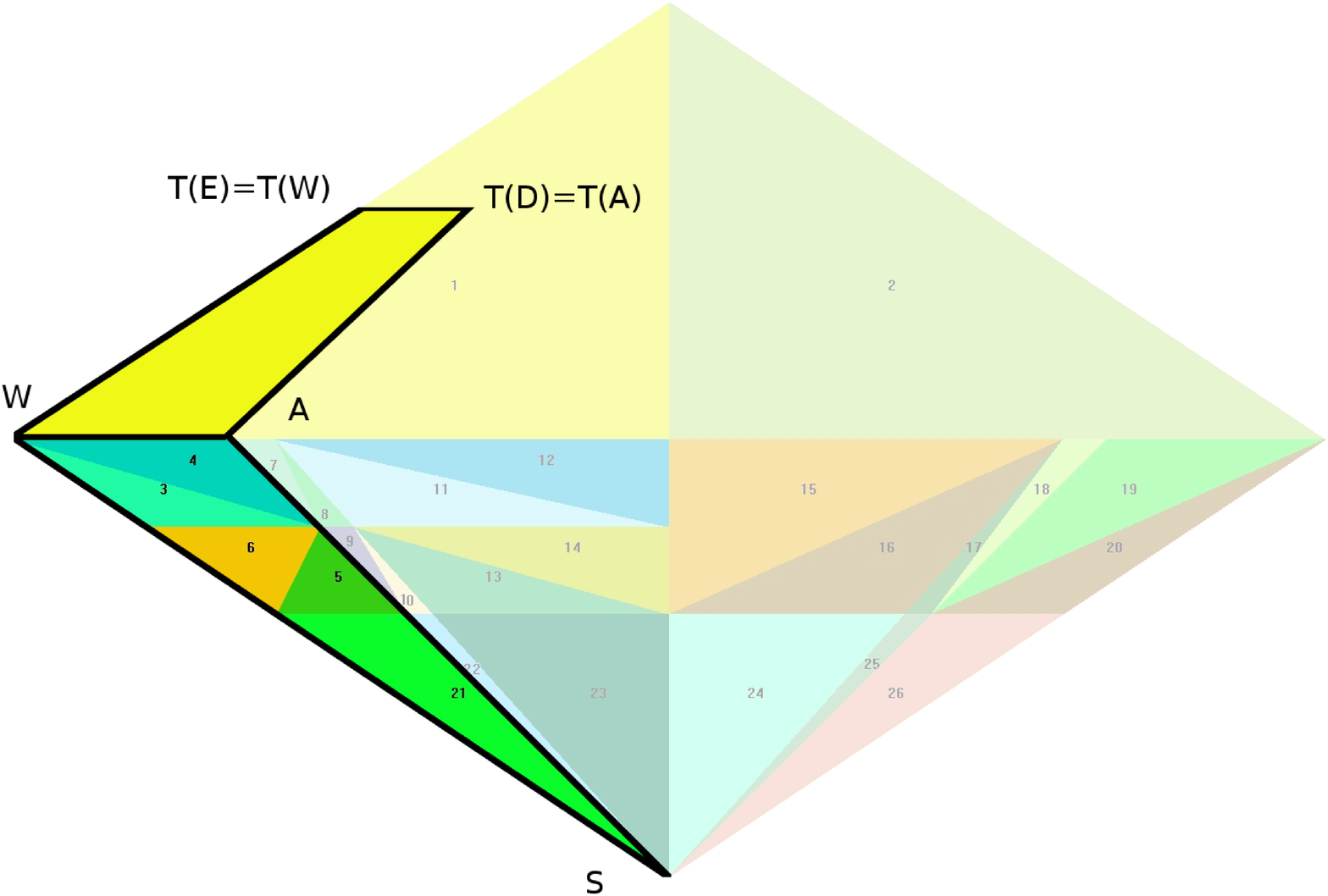}
\caption{The images of $(WAS)$ on the left and $(DES))$ on the right. The bold lines are the image of those in Fig. \ref{fig:left-right}.}\label{fig:image-left-right}
\end{figure}

Similarly one checks that $DES$ and $WAS$ (see Fig. \ref{fig:left-right}) are 
both mapped to subsets of $WAS \cup NEW$ (see Fig. \ref{fig:image-left-right}). This
establishes property (\ref{key-right}) and prepares the proof of (\ref{key-left}).

\subsection{Orbits not contained in $A^tB^tS\cup C^cD^cS$} By the above remarks,
such orbits must eventually land in $NEW$ --in which case they converge to the fixed point $N$; or
enter $WAS$ and stay there for ever.

Let us show that an orbit which is confined to $WAS$ converges to a fixed point.
$WAS$ is the union of $5$ triangles on each of which the map is affine: the triangles
numbered $3,4,5,6$ and $21$ on Fig. \ref{fig:partition}.

As pictured in Fig. \ref{fig:preimage-top}, the triangles $3$ and $4$ are mapped into $NEW$
so our orbit cannot enter them. As can be seen in Fig. \ref{fig:image-left-right}, the 
triangle $21$ (i.e., $W^cA^cS$) is mapped into itself. Moreover the segment $[W^cS]$ is 
made of fixed points, while the transverse direction
is contracted like $[A^cS]$. It follows that all points in this triangle eventually converge to
some fixed point in $[W^cS]$.

Triangle $6$, i.e., $W^tW^cA^t$, is mapped to $WW^cA$ which is contained in
$W^tW^cA^t\cup T^{-1}NEW$. But
 $$
    T'(x,y) = \left(\begin{matrix} \frac54 & \frac58\\ 0 & \frac53\end{matrix}\right)
 $$
with eigenvalues are $5/3$ and $5/4$: it is expanding. Thus all points in triangle 6, except the
fixed point $W^c$ are eventually mapped into $NEW$ under positive iteration.

We consider triangle $5$, i.e., $W^cA^cA^t$. Using Fig. \ref{fig:image-left-right}, 
observe that points that exit this triangle cannot re-enter it.
Hence it is enough to analyze orbits that stay in $W^cA^cA^t$ forever. However,
on that triangle, 
 $$
    T'(x,y) = \left(\begin{matrix} 2 & -\frac12\\ -1 & \frac32\end{matrix}\right)
 $$
with eigenvalues: $2.5$ and $1$, the latter with eigenvector $\left(\begin{matrix} 1\\ 2\end{matrix}\right)$. As $T(W^c)=W^c$, this gives a segment of fixed points for
the affine map and everything else is eventually mapped outside of the
triangle.

This completes the proof of property (\ref{key-left}) and therefore of all
the key properties.

\section{Proof of the Theorem}

In this section, we deduce the Main Theorem from the key properties (\ref{key-NS})-(\ref{key-left}).

First, let us note that the only aperiodic ergodic, invariant measures
are carried by the compact invariant set
 $$
    K := \bigcap_{n\geq0} T^{-n}(ABS\cup CDS).
 $$
Indeed, by key properties (\ref{key-half}), (\ref{key-folding}), (\ref{key-right}) and (\ref{key-right}) and the additional remark below them, all orbits which do not stay in $ABS\cup CDS$ for ever
converge to fixed points.

We set $K_0:=K\setminus\{S\}$ and denote the $2$-shift by $\sigma:\Sigma\to\Sigma$ with $\Sigma:=\{0,1\}^{\mathbb N}$. The key fact is:

\begin{lemma}
$(T,K_0)$ is an entropy-preserving extension of a subset of the $2$-shift according to
 $$
    \gamma:x\in K_0\longmapsto (\eps_n)_{n\geq0}\in\Sigma
      \text{ where }\eps_n=0\iff T^n(x)\in ABS.
 $$
More precisely,  $h(T,\mu)=h(\sigma,\gamma_*\mu)$ for any invariant probability measure $\mu$ of $(T,K_0)$.
\end{lemma}

\begin{proof}
$\gamma$ is clearly continuous and satisfies $\gamma\circ T=\sigma\circ\gamma$ on $K_0$.
We claim that $\gamma^{-1}\gamma(x)$ is a $C^1$ curve starting from $S$, containing $x$
and whose tangent is everywhere contained in $C^s$ (we call such curves \emph{vertical}).

Indeed, define, for any $n\geq0$, $C_n(x):=\{y\in Q:\forall 0\leq k\leq n\;
f^ny\in\gamma_k(x)\}$ where $\gamma_k(x)\in\{ABS,CDS\}$ is characterized by
$T^k(x)\in\gamma_k(x)$. Note that $T(ABS),T(CDS)\supset ABS\cup CDS$. 
Hence $T^nC_n(x)=ABS$ or $CDS$. 

According to (\ref{key-ABS})-(\ref{key-CDS}), $\mathcal C^s$ is stable by 
$T|ABS$ and $T|CDS$ and contains the vertical boundary lines of these two triangle (with equations $x=\kappa y$ with $|\kappa|\in\{0.9,1\}$, see Table \ref{tab:map-points}). 
Thus $C_n(x)$ is bounded by two
vertical curves and some graph $y=\phi(x)$. Moreover $(T^n|C_n(x))^{-1}$ is strongly contracting
horizontally by (\ref{key-transverse}). Therefore $\bigcap_{n\geq0} C_n(x)$ is a vertical curve
containing $x$.

To conclude, observe that this vertical curve is $\gamma^{-1}\gamma(x)$ and that the
restriction of $T^n$ to this set for any $n\geq0$ is a homeomorphism. Hence $h_\top(T, \gamma^{-1}\gamma(x))=0$.
It follows from a result of Bowen \cite{Bowen} that $\gamma$ preserves the entropy (the
requirement of compactness can be fulfilled by replacing $K_0$ with its image under
$(x,y)\mapsto(x/y,y)$, compactified by the addition of $\Sigma$).
\end{proof}

\begin{lemma}
$T$ satisfies: $h_\top(T)=\log 2$.
\end{lemma}

\begin{proof}
As a continuous map, $T$ satisfies:
$h_\top(T)=\sup_\mu h(T,\mu)$. The above lemma shows that one
can restrict this supremum to invariant measures carried by
$K_0$ and that these measures have entropy at most $\log 2$,
proving $h_\top(T)\leq\log 2$. For the reverse implication,
we use that $T|A^cB^cS$ and $T|C^cD^cS$ are linear,
multiplying the $y$-coordinate by $1/2$ and $2$ respectively.
It follows that if $(x_k,y_k):=T^k(x_0,y_0)$ belongs to these
two small triangles near $S$ for $0\leq k<n$, then:
 $$
    \log y_n = \log y_0 + \log 2 \cdot \sum_{k=0}^{n-1}\operatorname{sign}(x_k).
 $$
Let $M$ be a large integer. Let
 $$
   \Sigma_M:=\biggl\{\alpha\in\Sigma:\forall p<q\;
      \bigl|\sum_{k=p}^q(\alpha_k-\frac12)\bigr| \leq M\biggr\}.
 $$
It is clearly compact and invariant, i.e., a subshift.  We claim that:
$\lim_{M\to\infty} h_\top(\Sigma_M)=\log 2$. 

Indeed, let $B(N):=\{A\in\{0,1\}^{2N}:\sum_{k=0}^{2N-1}A_k=N\}$ and
 $$
    \Sigma'_M:=\bigcup_{k=0}^{2N-1}\{A^1A^2\dots:A^n\in B(N)\}
  $$
This is a subshift of $\Sigma_{2N}$ with entropy $\log\#B(N)/2N$ which
converges to $\log 2$. The claim is proven.

Let
 $$
    X_M:=\{(\alpha,s)\in\Sigma_M\times\{-M,\dots,M\}:\forall p\in\NN\;
       \bigl|\frac s2+\sum_{k=0}^p(\alpha_k-1/2)\bigr|\leq \frac M2\}
 $$
and define $F_M:X_M\to X_M$ by $F_M(\alpha,s)=(\sigma(\alpha),s+(\alpha_0-1/2))$.
It is easy to check that this is a well-defined, finite, topological extension of $\Sigma_M$.
Also it can be embedded into $K_0\cap\{2^{-M}y^c\leq y\leq y^c\}$ by the
map $\iota_M$ defined by:
 $$
    \iota_M(\alpha,s) = (y(s)x(\alpha),y(s))
 $$
where
 $$
    x(\alpha):=-\frac{19}{20}\sum_{n\geq0} \frac{\sigma_0\dots\sigma_n}{20^n}
      \text{ and }y(s):=y^c2^{s-M/2}
 $$
Embedding a measure maximizing entropy of $F_M$ into $T$ through $\iota_M$ and
letting $M\to\infty$ shows that $h_\top(T)=\log2$.
\end{proof}

\begin{proposition}
$T$ has no invariant probability measure with entropy $\log 2$.
\end{proposition}

\begin{proof}
We proceed by contradiction, assuming the existence of such a measure,
say $\mu$. By the previous lemma, it is supported by $K_0$. Hence
$\gamma_*\mu$ is an invariant probability measure of the full shift
with entropy $\log 2$. It must be the $(\frac12,\frac12)$-Bernoulli
measure. Thus, noting as above $(x_k,y_k):=T^k(x_0,y_0)$, we have:
 $$
    \mu-\forall (x_0,y_0)\qquad 
      \sup_{n\geq0} \sum_{k=0}^n \operatorname{sign}(x_k)=\infty.
 $$
Now,
 $$
    \log y_n = \log y_0 + \sum_{k=0}^{n-1} \log\frac{y_{k+1}}{y_k}
 $$
and $y_{k+1}\geq 2^{\operatorname{sign}(x_k)}y_k$ by key properties (\ref{key-ABS})-(\ref{key-CDS}).
Hence, $\sup_{n\geq0} y_n=\infty$, a contradiction.
\end{proof}


\begin{thebibliography}{99999}
\bibitem{Bowen}
R. Bowen, Entropy for group endomorphisms and homogeneous spaces, {\it Trans. Amer. Math. Soc.}  {\bf 153}  (1971), 401--414.

\bibitem{BuzziSIM}
J. Buzzi, Intrinsic ergodicity of smooth interval maps, {\it  Israel J. Math.} {\bf  100}  (1997), 125--161.

\bibitem{BuzziAffine}
J. Buzzi, Intrinsic ergodicity of affine maps in $[0,1]^d$, {\it  Monatsh. Math.} {\bf 124}  (1997),  no. 2, 97--118.

\bibitem{BuzziPWAH}
J. Buzzi, Measures of maximal entropy for piecewise affine surface homeomorphisms,
{\it Ergod. th. dynam. syst.}, (to appear).

\bibitem{Hofbauer}
F. Hofbauer, On intrinsic ergodicity of piecewise monotonic transformations with positive entropy, {\it Israel J. Math.} {\bf 34}  (1979), 213--237; II. {\it Israel J. Math.} {\bf 38}  (1981), 107--115. 

\bibitem{KruglikovRypdal1}
B. Kruglikov, M. Rypdal, A piece-wise affine contracting map with positive entropy,
{\it  Discrete Contin. Dyn. Syst.} {\bf 16}  (2006),  393--394.

\bibitem{KruglikovRypdal2}
B. Kruglikov, M. Rypdal, Entropy via multiplicity, {\it  Discrete Contin. Dyn. Syst.} {\bf 16}  (2006),  no. 2, 395--410.

\bibitem{Newhouse}
S. Newhouse, Entropy and volume, {\it Ergodic Theory Dynam. Systems} {\bf 8$\sp *$}  (1988),  Charles Conley Memorial Issue, 283--299.

\bibitem{NewhousePerso}
S. Newhouse, {\it personal communication,} 1999.

\bibitem{RuetteEx}
S. Ruette, Mixing $C\sp r$ maps of the interval without maximal measure, {\it Israel J. Math.} {\bf 127}  (2002), 253--277.

\bibitem{Yomdin}
Y. Yomdin, Volume growth and entropy, {\it Israel J. Math.} {\bf 57}  (1987),  285--300.

\end{thebibliography}
\end{document}